\documentclass{amsart}
\usepackage{amscd,amssymb,latexsym,epsfig,amsmath,xypic,supertabular}
\begin{document}
\newcommand{\mmbox}[1]{\mbox{${#1}$}}
\newcommand{\proj}[1]{\mmbox{{\mathbb P}^{#1}}}
\newcommand{\affine}[1]{\mmbox{{\mathbb A}^{#1}}}
\newcommand{\Ann}[1]{\mmbox{{\rm Ann}({#1})}}
\newcommand{\caps}[3]{\mmbox{{#1}_{#2} \cap \ldots \cap {#1}_{#3}}}
\newcommand{\N}{{\mathbb N}}
\newcommand{\Z}{{\mathbb Z}}
\newcommand{\R}{{\mathbb R}}
\renewcommand{\k}{\Bbbk}
\newcommand{\p}{{\mathbb P}}
\newcommand{\A}{{\mathcal A}}
\newcommand{\CC}{{\mathcal C}}
\newcommand{\C}{{\mathbb C}}

\newcommand{\Tor}{\mathop{\rm Tor}\nolimits}
\newcommand{\Ext}{\mathop{\rm Ext}\nolimits}
\newcommand{\Hom}{\mathop{\rm Hom}\nolimits}
\newcommand{\im}{\mathop{\rm Im}\nolimits}
\newcommand{\rank}{\mathop{\rm rank}\nolimits}
\newcommand{\supp}{\mathop{\rm supp}\nolimits}
\newcommand{\coker}{\mathop{\rm coker}\nolimits}
\newtheorem{defn0}{Definition}[section]
\newtheorem{prop0}[defn0]{Proposition}
\newtheorem{conj0}[defn0]{Conjecture}
\newtheorem{thm0}[defn0]{Theorem}
\newtheorem{lem0}[defn0]{Lemma}
\newtheorem{corollary0}[defn0]{Corollary}
\newtheorem{example0}[defn0]{Example}
\newtheorem{remark0}[defn0]{Remark}
\newtheorem{que0}[defn0]{Question}

\newenvironment{defn}{\begin{defn0}}{\end{defn0}}
\newenvironment{prop}{\begin{prop0}}{\end{prop0}}
\newenvironment{conj}{\begin{conj0}}{\end{conj0}}
\newenvironment{thm}{\begin{thm0}}{\end{thm0}}
\newenvironment{lem}{\begin{lem0}}{\end{lem0}}
\newenvironment{que}{\begin{que0}}{\end{que0}}
\newenvironment{cor}{\begin{corollary0}}{\end{corollary0}}
\newenvironment{exm}{\begin{example0}\rm}{\end{example0}}
\newenvironment{remark}{\begin{remark0}\rm}{\end{remark0}}

\newcommand{\defref}[1]{Definition~\ref{#1}}
\newcommand{\propref}[1]{Proposition~\ref{#1}}
\newcommand{\thmref}[1]{Theorem~\ref{#1}}
\newcommand{\lemref}[1]{Lemma~\ref{#1}}
\newcommand{\corref}[1]{Corollary~\ref{#1}}
\newcommand{\exref}[1]{Example~\ref{#1}}
\newcommand{\secref}[1]{Section~\ref{#1}}
\newcommand{\rmkref}[1]{Remark~\ref{#1}}
\newcommand{\queref}[1]{Question~\ref{#1}}
\title {High rank linear syzygies on low rank quadrics}

\author{Hal Schenck}
\thanks{Schenck supported by NSF 0707667, NSA H98230-07-1-0052}\address{Schenck: Mathematics Department \\ University of Illinois\\
  Urbana \\ IL 61801\\ USA}
\email{schenck@math.uiuc.edu}

\author{Mike Stillman}
\thanks{Stillman supported by NSF 0810909}
\address{Stillman: Mathematics Department \\ Cornell University \\
  Ithaca \\ NY 14850\\ USA}
\email{mike@math.cornell.edu}


\begin{abstract}
\noindent 
We study the linear syzygies of a homogeneous ideal
$I \subseteq S = Sym_{\k}(V)$, focussing on the graded betti numbers
\[
b_{i,i+1} = \dim_{\k}Tor_i(S/I, \k)_{i+1}.
\]
For a variety $X$ and divisor $D$ with $V = H^0(D)$, what 
conditions on $D$ ensure that $b_{i,i+1} \ne 0$? In \cite{e3},
Eisenbud shows that a decomposition $D\! \sim\! A\!+\!B$ such that $A$ and $B$ have 
at least two sections gives rise to determinantal equations 
(and corresponding syzygies) in $I_X$; and in \cite{e1} conjectures
that if $I_2$ is generated by quadrics of rank $\le 4$, then 
the last nonvanishing $b_{i,i+1}$ is a consequence of such
equations. We describe obstructions to this 
conjecture and prove a variant. The obstructions arise
from toric specializations of the Rees algebra of Koszul cycles,
and we give an explicit construction of toric varieties with 
minimal linear syzygies of arbitrarily high rank. This gives
one answer to a question posed by Eisenbud and Koh in \cite{ek}. 
\end{abstract}
\maketitle
\vspace{-.2in}
\section{Introduction}\label{sec:intro}
Let $I$ be a homogeneous ideal in $S= Sym_{\k}(V)$, with $\k$
algebraically closed and characteristic zero; we are 
primarily interested in the case that $I$ 
is the ideal of an irreducible, nondegenerate variety in $\p(V)$, 
and in the graded betti numbers 
\[
b_{ij} = dim_{\k}Tor_i(S/I,\k)_j.
\]
\begin{defn}
The length of the 2-linear strand $2LP(S/I) = \max\{i | b_{i,i+1} \ne 0\}.$
\end{defn}
\noindent In particular, $I$ has a $2$-linear $n^{th}$--syzygy iff $2LP(S/I)\ge n+1$.
\begin{exm}\label{exm:twistedcubic}
The twisted cubic $I_{C} \subseteq S=\k[x,y,z,w]$ has resolution

\begin{small}
\[
0 \longrightarrow S(-3)^2 \xrightarrow{\left[ \!
\begin{array}{cc}
-z & w\\
y & -z \\
-x & y
\end{array}\! \right]} S(-2)^3
\xrightarrow{\left[ \!\begin{array}{ccc}
y^2-xz& yz-xw& z^2-yw
\end{array}\! \right]}
 S \longrightarrow S/I_{C}
\]
\end{small}
So $2LP(S/I_{C}) = 2$; in betti diagram notation \cite{e3} the $b_{ij}$ for $S/I_X$ are:
\begin{small}
$$
\vbox{\offinterlineskip 
\halign{\strut\hfil# \ \vrule\quad&# \ &# \ &# \ &# \ &# \ &# \
&# \ &# \ &# \ &# \ &# \ &# \ &# \
\cr
total&1&3&2\cr
\noalign {\hrule}
0&1 &--&--&\cr
1&--&3 &2 &\cr
\noalign{\bigskip}
\noalign{\smallskip}
}}
$$
\end{small}
\end{exm}
The Eagon-Northcott complex gives a free resolution for $I_{C}$, 
which is determined by the maximal minors of the matrix of
first syzygies; this simple example provides the key intuition: 
if $D\sim A+B$, then $h \in H^0(A)$ and $g \in H^0(B)$ yield 
an element $f = g\cdot h \in H^0(D)$. Factoring a divisor of degree
three on $\p^1$ into divisors of degree one and two yields a $2 \times
3$ matrix for $(s,t) \otimes (s^2,st,t^2)$:
\[
\left[ \!
\begin{array}{cc}
st^2 & t^3\\
s^2t & st^2 \\
s^3 & s^2t
\end{array}\! \right]
\]
So $I_{C}$ contains the $2 \times 2$ minors of a $2 \times
3$ matrix of linear forms. This matrix is of a special type:
\begin{defn}
A matrix of linear forms is $1-generic$ if it has no zero entry, and
cannot be transformed by row and column operations to have a zero entry.
\end{defn}
In \cite{e3}, Eisenbud shows that for a reduced irreducible
nondegenerate linearly normal curve $C \subseteq \p^r$, there is a 1-generic
$p \times q$ matrix of linear forms whose $2 \times 2$ minors vanish
on $C$ iff there exist line bundles ${\mathcal L}_1$ and ${\mathcal
  L}_2$ such that
\[
{\mathcal O}_C(1) \simeq {\mathcal L}_1 \otimes {\mathcal L}_2,
\]
with $h^0({\mathcal L}_1) \ge p$ and $h^0({\mathcal L}_2) \ge q$.
Combining this with a result of Eisenbud 
that the ideal of $2 \times 2$ minors of 
a 1-generic $p \times q$ matrix has $2LP \ge p+q-3$ 
(see \cite{ks}) leads to:

\begin{conj}\label{conj:GG}$[$Eisenbud, \cite{e1}$]$
Let $\k$ be algebraically closed of characteristic zero,
and $I \subseteq S$ be a prime ideal containing no linear form,
such that $I_2$ is spanned by quadrics of rank at most four. 
If $2LP(S/I)=n$, then $I$ contains the $2 \times 2$ minors of a
1-generic $p \times q$ matrix, with $p+q-3 = n$.
\end{conj}
With the additional hypotheses
that $S/I$ is normal, Gorenstein, of dimension $2$ and degree $2r$,
Conjecture~\ref{conj:GG} specializes to Green's conjecture \cite{g1}. 
Our first result (Theorem~\ref{thm:first} below) provides an 
infinite class of counterexamples to Conjecture~\ref{conj:GG}.
To state the result we need:
\begin{defn}
For a homogeneous ideal $I \subseteq S$, let $F_\bullet$ 
denote the subcomplex of the minimal free resolution of $S/I$
\[
0 \longleftarrow S/I_2 \longleftarrow S  \longleftarrow F_1 \otimes
S(-2) \longleftarrow F_2 \otimes S(-3)  \longleftarrow F_3 \otimes
S(-4) \longleftarrow \cdots
\]
Let $f$ be a $2$-linear $n^{th}$ syzygy of $I$ (henceforth, ``linear $n^{th}$ syzygy''). The rank of $f$ is the 
dimension of the smallest vector space $G$ such that the 
diagram below commutes:
\[
\xymatrix{
F_{n} \otimes S(-n-1)  & F_{n+1} \otimes S(-n-2)\ar[l]\\
G \otimes S(-n-1) \ar[u] &  f \otimes S(-n-2)\ar[l]\ar[u] }
\]
\end{defn}
\vskip .1in
\begin{thm}\label{thm:first}
For any odd $n$, there exists an arithmetically Cohen-Macaulay 
toric ideal generated by $n$ quadrics of rank $\le 4$, 
with only one linear first syzygy, of rank $n$.
For any even $n$, there exists an arithmetically Gorenstein toric ideal, 
generated by $n$ quadrics of rank $\le 4$, 
and a $\frac{n}{2}$-ic, with only one linear first syzygy, of rank $n$.
\end{thm}
These ideals are counterexamples to Conjecture~\ref{conj:GG}: 
$2LP(S/I)=2$ should force $b_{23} \ge 2$. 
Roughly speaking, the problem is that the condition that $I_2$ is 
generated by quadrics of rank at most four does not guarantee
that there cannot be $2$-linear syzygies (possibly all) of comparatively
high rank. Hence, some additional hypothesis, such as requiring
a top linear syzygy of low rank, is necessary. 

\begin{thm}\label{thm:second}
Let $\k$ be algebraically closed of characteristic zero,
and $I \subseteq S$ be a prime ideal containing no linear form,
such that $2LP(S/I)=n+1$. 
\begin{enumerate}
\item If $I$ has a linear $n^{th}$ syzygy of rank $n+2$, then $I$
  contains the $2 \times 2$ minors of a 1-generic $2 \times (n+2)$ matrix. 
\item If $I$ has a linear $n^{th}$ syzygy of rank $n+3$, then $I$
  contains the $4 \times 4$ Pfaffians of a skew-symmetric 
1-generic $(n+4) \times (n+4)$ matrix.
\item If $I$ has a linear $n^{th}$ syzygy of rank $n+3$ {\em and} is 
a semigroup ideal, then $I$ contains the $2 \times 2$ minors of a 
1-generic $p \times (n-p+4)$ matrix.
\item If $I$ has no linear $n^{th}$ syzygies of rank $\le n+3$, then 
$I$ does not contain the $2 \times 2$ minors of a 1-generic $p \times q$ matrix, with $p+q-3 = n+1$.
\end{enumerate}
\end{thm}
Note that if $I$ has a linear $n^{th}$ syzygy of rank $\le n+1$, then
$I$ cannot be prime.
Our main tool in studying Conjecture~\ref{conj:GG} is Koszul 
homology. Recall that $\Tor_m(S/I, k)_{m+1}$ may be computed
as the homology of
\[
\bigwedge^{m+1} V \stackrel{\partial}{\rightarrow} \bigwedge^m V \otimes V
\stackrel{\partial}{\rightarrow} \bigwedge^{m-1} V \otimes (S/I)_2
\]
where $V \simeq (S/I)_1 = S_1$. Without the rank conditions, the first
two parts of Theorem~\ref{thm:second} appear in \cite{ks}.

Perhaps our most interesting result arises from the fact that the ideals
which appear in Theorem~\ref{thm:first} are toric specializations
of the Rees algebra of Koszul cycles. Such Rees algebras have been
previously studied in \cite{htz}, \cite{hu}, \cite{k86} and \cite{w}. 
In \cite{ek}, Eisenbud and Koh ask ``Under what conditions does
a module with a linear $k^{th}$ syzygy specialize to one with a linear
$k^{th}$ syzygy?'' We give one answer to this question: in \S 4 we prove:
\begin{thm}\label{thm:third}
Let $\Delta$ be an oriented $n$--dimensional 
pseudomanifold on $d$--vertices, with top homology class
$H_n(\Delta)$. Then there exists a toric variety $X=V(I)$ of 
dimension $d$, such that $I$ has a {\em minimal} 
linear $n^{th}$--syzygy of rank $d$. If ${\bf m}$ is
the multidegree of this syzygy, then the complex 
$\Delta_{\bf m}$ which computes $Tor_n^S(I,\k)_{\bf m}$
is homotopic to $\Delta$.
\end{thm}
It seems reasonable to ask if Conjecture~\ref{conj:GG} holds with
additional geometric constraints, and we explore this in \S 5. 
For example, taking linear sections of the 
varieties appearing in Theorem~\ref{thm:first} yields smooth, 
projectively normal curves for which Conjecture~\ref{conj:GG} fails,
but the associated divisor is special.  
\section{Linear first syzygies}
We start by analyzing the linear first syzygies, computed as the homology of:
\[
\bigwedge^{3} V \stackrel{\partial_2}{\rightarrow} \bigwedge^2 V \otimes V
\stackrel{\partial_1}{\rightarrow} V \otimes (S/I)_2
\]
where (unless otherwise noted) $V \simeq S_1$, with a basis element 
$e_i \in V$ mapping to $x_i \in S$. The differential $\partial_1$ is defined via
\[
\partial_1(e_i \wedge e_j \otimes y_{\widehat{ij}}) = y_{\widehat{ij}}\cdot (x_i e_j
-x_je_i),
\]
where $y_{\widehat{ij}}$ is an indeterminate linear form in $S$ (the reason for the $\mbox{ } \widehat{}\mbox{ } $ will appear in \S 4). 
A particular class in Koszul homology representing a rank $d$ syzygy will be 
supported on a subspace $W$ of dimension $d$:
\[
\bigwedge^2 W \otimes V \subseteq \bigwedge^2 V \otimes V.
\]
\begin{exm}\label{exm:3vars}
Suppose $\k^3 \simeq W \subseteq V$, so that an element of $\bigwedge^2 W \otimes
V$ may be written as $
\omega = e_1 \wedge e_2 y_{\widehat{12}}-  e_1 \wedge e_3 y_{\widehat{13}}
+ e_2 \wedge e_3 y_{\widehat{23}}.$
\begin{figure}[h]
\epsfig{file=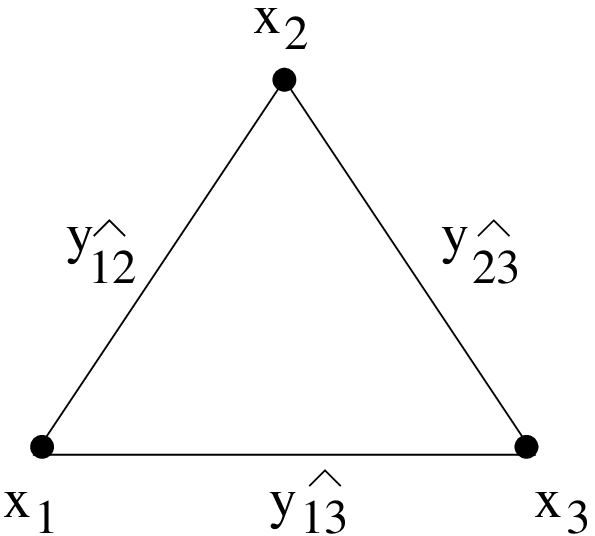,height=1.2in,width=1.2in}
\end{figure}
\vskip -.05in
\noindent Thus, 
\[
\partial(\omega) = y_{\widehat{12}} (x_1e_2-x_2e_1)- y_{\widehat{13}}(x_1e_3 - x_3e_1)
+ y_{\widehat{23}}(x_2e_3 - x_3e_2)
\]
is a cycle for $S/I$ iff 
$I' = \langle y_{\widehat{13}}x_3-y_{\widehat{12}} x_2, y_{\widehat{12}} x_1 - y_{\widehat{23}}x_3, y_{\widehat{23}}x_2 - y_{\widehat{13}} x_1\rangle \subseteq I$.
Now, $I'$ is generated by the $2 \times 2$ minors of 
\[
\left[ \!
\begin{array}{ccc}
x_1 & x_2 & x_3\\
y_{\widehat{23}} &  y_{\widehat{13}} & y_{\widehat{12}}
\end{array}\! \right].
\]
If $V(I')$ is irreducible and nondegenerate (note that some of
the $y_{\widehat{ij}}$ could vanish), then $2LP(S/I') =2$, so Conjecture \ref{conj:GG} 
holds for $S/I'$.
As observed by Schreyer (\cite{s},
Lemma 4.3) this process generalizes; let 
\begin{equation}\label{schreyerEq}
\omega = \sum\limits_{1\le i < j \le n} (-1)^{i+j+1}y_{\widehat{ij}}e_i \wedge e_j.
\end{equation}
\end{exm}
\noindent When $d=4$, this yields Pfaffians, but for $d \ge 5$, 
interesting behavior occurs. The term $(-1)^{i+j+1}$ simplifies the connection to the Rees algebra in \S 4.
\begin{exm}\label{exm:5vars}
If $d = 5$, the quadratic conditions necessary for 
$\partial(\omega)$ to represent a class in homology have rank
eight. 
\[
\left[ \!
\begin{array}{cccccc}
0 & -y_{\widehat{12}} & y_{\widehat{13}} & -y_{\widehat{14}} &y_{\widehat{15}}\\
y_{\widehat{12}} &0 &-y_{\widehat{23}} & y_{\widehat{24}} & -y_{\widehat{25}}\\
-y_{\widehat{13}} &y_{\widehat{23}}&0 & -y_{\widehat{34}} & y_{\widehat{35}}\\
y_{\widehat{14}} &-y_{\widehat{24}}&y_{\widehat{34}} &0 &-y_{\widehat{45}}\\
-y_{\widehat{15}} &y_{\widehat{25}}& -y_{\widehat{35}} &y_{\widehat{45}}&0
\end{array}\! \right]
\cdot
\left[ \!
\begin{array}{c}
x_1\\
x_2\\
x_3\\
x_4\\
x_5
\end{array}\! \right]=0
\]
\vskip .05in
\noindent Specializing to a simple cycle 
by setting $y_{\widehat{13}} = y_{\widehat{14}} = y_{\widehat{24}} = y_{\widehat{25}} = y_{\widehat{35}} = 0$
\begin{figure}[h]
\epsfig{file=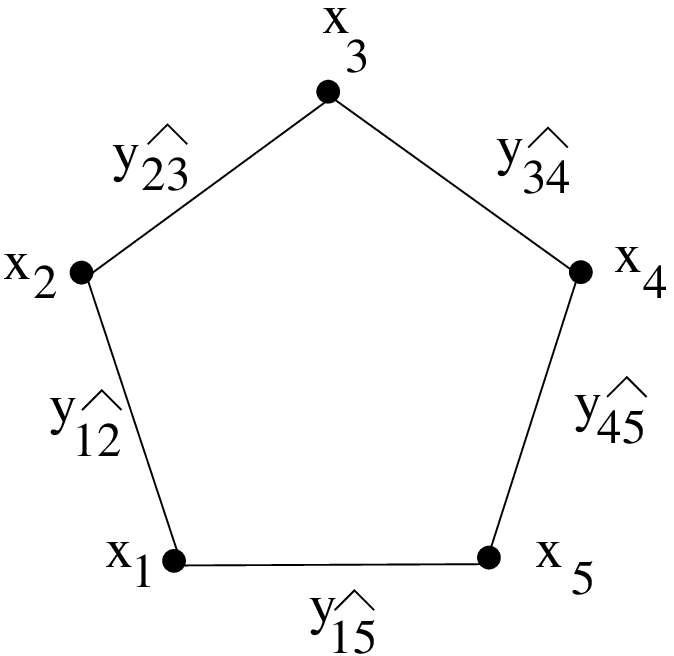,height=1.4in,width=1.4in}
\end{figure}
\vskip .05in
yields an ideal generated by quadrics of rank at most four:
\[
\left[ \!
\begin{array}{cccccc}
0 & -y_{\widehat{12}} & 0 & 0 &y_{\widehat{15}}\\
y_{\widehat{12}} &0 &-y_{\widehat{23}} & 0 & 0 \\
0 &y_{\widehat{23}}&0 & -y_{\widehat{34}} & 0\\
0 & 0& y_{\widehat{34}} &0 &-y_{\widehat{45}}\\
-y_{\widehat{15}} & 0& 0 & y_{\widehat{45}}&0
\end{array}\! \right]
\cdot
\left[ \!
\begin{array}{c}
x_1\\
x_2\\
x_3\\
x_4\\
x_5
\end{array}\! \right]=0
\]
In contrast to the $d=4$ case, the resulting ideal is prime, and 
defines a (singular) 5-fold $X \subseteq \p^9$, with graded betti numbers:
\begin{small}
$$
\vbox{\offinterlineskip 
\halign{\strut\hfil# \ \vrule\quad&# \ &# \ &# \ &# \ &# \ &# \
&# \ &# \ &# \ &# \ &# \ &# \ &# \
\cr
total&1&5&12&10&2\cr
\noalign {\hrule}
0&1 &--&--& --&--&    \cr
1&--&5 &1 & --&--&        \cr
2&--&-- &11 & 10&1&        \cr
3&--&-- &-- &--&1&         \cr
\noalign{\bigskip}
\noalign{\smallskip}
}}
$$
\end{small}
Since $2LP =2$, Conjecture \ref{conj:GG} would require that
$I_X$ contain the $2\times 2$ minors of a $1-$generic $2 \times 3$ 
matrix. Such a matrix would force existence of a pair of linear first syzygies,
so this is a counterexample to Conjecture~\ref{conj:GG}. The linear syzygy
we obtain is by construction of rank five, so is of rank 
too high to be a linear first syzygy for an ideal of $2\times 2$
minors. 

$X$ is arithmetically Cohen-Macaulay with singularities in codimension $3$;
quotienting by four generic linear forms yields a smooth,
projectively normal curve of genus seven in $\p^5$, generated by
quadrics of rank $\le 4$. This curve is a projection of a general 
canonical curve $C \subseteq \p^6$ from a
point $ p \in C$; $h^1(K_C-p) = 1$ so $K_C-p$ is special. We analyze this
in more detail in \S 5.
\end{exm}

\begin{exm}\label{exm:6vars}
The last explicit example we give is for the $d = 6$ case;
specializing to a cycle yields the ideal $I_6$:
\[
\left[ \!
\begin{array}{cccccc}
0 & -y_{\widehat{12}} & 0 & 0 & 0 &-y_{\widehat{16}}\\
y_{\widehat{12}} &0 &-y_{\widehat{23}} & 0 & 0 &0\\
0 &y_{\widehat{23}}&0 & -y_{\widehat{34}} & 0&0\\
0 & 0& y_{\widehat{34}} &0 &-y_{\widehat{45}}&0\\
0 & 0& 0       &y_{\widehat{45}} &0 &-y_{\widehat{56}}\\
y_{\widehat{16}} & 0& 0 &0& y_{\widehat{56}}&0
\end{array}\! \right]
\cdot
\left[ \!
\begin{array}{c}
x_1\\
x_2\\
x_3\\
x_4\\
x_5\\
x_6
\end{array}\! \right]=0
\]
As with $I_4$, the ideal $I_6$ is not prime; 
$I_6 = \langle x_1,\ldots,x_6\rangle \cap J_6$, with 
$J_6$ prime of codimension 5. 
\[
J_6 = I_6 + \langle y_{\widehat{12}}y_{\widehat{34}}y_{\widehat{56}}+y_{\widehat{23}}y_{\widehat{45}}y_{\widehat{16}}\rangle,
\]
Taking hyperplane sections yields a smooth curve 
$X \subseteq \p^6$ of degree $21$ and genus $21$ with 
betti numbers identical to those of $J_6$; in
particular $X$ is Gorenstein:
\begin{small}
$$
\vbox{\offinterlineskip 
\halign{\strut\hfil# \ \vrule\quad&# \ &# \ &# \ &# \ &# \ &# \
&# \ &# \ &# \ &# \ &# \ &# \ &# \
\cr
total&1&7&22&22&7&1\cr
\noalign {\hrule}
0&1 &--&--& --&--& --&   \cr
1&--&6 &1 & --&--&--&        \cr
2&--&1 &21 & 21&1&--&        \cr
3&--&-- &-- &1&6&--&         \cr
4&--&-- &-- &--&--&1&         \cr
\noalign{\bigskip}
\noalign{\smallskip}
}}
$$
\end{small}
\end{exm}
To analyze the situation for general $d$, let $Y$ be a generic $d \times d$ 
skew symmetric matrix with entries
$y_{\widehat{ij}}$ and $X$ a generic column vector with entries $x_i$. Write 
$I_d$ for the ideal generated by the entries of $Y \cdot X$, and let
$J_d = I_d + \mbox{Pfaff}(\mbox{det}(Y))$; when $d$ is odd $J_d=I_d$.
When $d$ is even, $S/J_d$ is Gorenstein of deviation two, and was
first studied by Huneke and Ulrich in \cite{hu}. In \cite{k86}, 
Kustin determined the minimal free resolution for both $I_d$ and 
$J_d$; in particular both are arithmetically Cohen-Macaulay. The ideals which
appear in Theorem~\ref{thm:first} are obtained from 
$I_d$ and $J_d$ by specializing all entries of $Y$ above the 
diagonal to zero, except the supradiagonal entries, and the 
entry in position $(1,d)$ (the ``corner'' entry). Write $\mathcal{I}_d$
and $\mathcal{J}_d$ for the specializations, and let $y_i$ 
denote the entry of the specialization in position $(i,i+1)$ and 
$y_d$ denote the ``corner'' entry.
\begin{lem}\label{lexlem}
The specializations described above are of codimension $d-1$.
\end{lem}
\begin{proof}
In lex order on $R=\k[y_1,\ldots,y_d,x_1,\ldots,x_d]$, the lead
monomials of the ideal generated by $(Y \setminus \mbox{row}_d(Y))\cdot
X$ are relatively prime, so the corresponding elements of $I_d$
generate a complete intersection of codimension $d-1$. 
As noted earlier, Kustin's results show that $I_d$ and $J_d$ are
Cohen-Macaulay of codimension $d-1$, so the specialization chosen
corresponds to a regular sequence (if it were not, the codimension 
would be less than $d-1$) which implies that the codimensions 
of $\mathcal{I}_d$ and $\mathcal{J}_d$ are also $d-1$.
\end{proof}
Theorem~\ref{thm:first} follows from Lemma~\ref{lexlem}. Since 
the ideals $\mathcal{I}_d$ and $\mathcal{J}_d$ are obtained from $I_d$ and
$J_d$ by quotienting with a regular sequence, the graded betti 
numbers of $\mathcal{I}_d$ and $\mathcal{J}_d$ are
identical to those appearing in the Kustin resolutions. 
In particular, $b_{23}(\mathcal{I}_k)= b_{23}(\mathcal{J}_d)= 1$. After
a linear change of coordinates $y_{1d}\mapsto -y_{1d}$ if $d$ is even, 
$\mathcal{I}_d$ and $\mathcal{J}_d$ are generated by binomials. Letting
 ${\bf e}_i$ denote the $i^{th}$ standard basis vector in $\k^{d+1}$, 
the parameterization $x_i = {\bf e}_i$, 
$y_{\widehat{ij}}={\bf e}_0-{\bf e}_i - {\bf e}_j$ shows they are toric.

\section{Proof of Theorem~\ref{thm:second}}
\noindent 
For an $n^{th}$ linear syzygy supported on a $d$-dimensional subspace $W \subseteq V = S_1$, let
$B$ be a basis for $W$, and for $I \subseteq B$ let 
\[
x_{\widehat{I}} = \bigwedge\limits_{j \in B \setminus I} x_j.
\]
To prove part $(1)$, suppose $I$ has an $n^{th}$ linear syzygy of rank $d=n+2$, 
supported on $W$, and let $S = \{x_1,\ldots,x_{n+2} \}$. Then 
\[
x_{\widehat{i}} = x_1 \wedge \cdots \wedge x_{i-1} \wedge x_{i+1} \wedge \cdots \wedge x_{n+2} \in
\bigwedge^{n+1}W.
\]
The syzygy is represented by 
\[
\omega = \sum\limits_{i=1}^{n+2}(-1)^{i+1} x_{\widehat{i}} \otimes y_i  \in
\bigwedge^{n+1}W \bigotimes V,
\]
where $y_i$ is an indeterminate linear form; since 
$i = [n+2] \setminus \{[n+2] \setminus i\}$, this agrees with the notation of
the previous section.
\[
\partial((-1)^{i+1}x_{\widehat{i}} \otimes y_i)=(-1)^{i+1}y_i
\big[ \sum\limits_{j=1}^{i-1}(-1)^{j-1}x_j x_{\widehat{ij}} +
\sum\limits_{j=i+1}^{n+2}(-1)^{j}x_jx_{\widehat{ij}} \big],
\]
it follows that $\partial(\omega)$ is a class in homology
iff the coefficients
of all the $x_{\widehat{ij}}$ are in $I_2$, which occurs exactly
when 
\[
(-1)^{i+1}y_i(-1)^{j}x_j + (-1)^{j+1}y_j(-1)^{i-1}x_i = 0.
\]
Hence, the $2 \times 2$ minors of $\phi=\left[ \!
\begin{array}{ccc}
x_1 & \cdots & x_{n+2}\\
y_1 & \cdots & y_{n+2}
\end{array}\! \right]
$
are in $I_2$, and $\omega$ represents a nontrivial class in homology iff 
some $2 \times 2$ minor is nonzero. The assumption
that $I$ is prime and nondegenerate implies that $\phi$ is one--generic. 
\vskip .1in
\noindent For part $(2)$, if $I$ has an $n^{th}$ linear syzygy of rank $n+3$, 
then let $W$ be a vector space of dimension $n+3$ with $B= \{x_1,\ldots,x_{n+3} \}$. 
The syzygy is represented by 
\[
\omega = \sum\limits_{1\le i < j \le n+3}(-1)^{i+j} x_{\widehat{ij}} \otimes y_{ij} \in
\bigwedge^{n+1}W \bigotimes V.
\]
A computation as above shows that $\partial(\omega)$ is a class in homology
when the coefficients of the $x_{\widehat{ijk}}$ are in $I_2$, which
happens if 
\[
y_{ij}x_k - y_{ik}x_j + y_{jk}x_i \in I_2, 
\]
so that $I_2$ contains the $4 \times 4$ Pfaffians of
\[
\left[ \!
\begin{array}{cccccc}
0       &x_1     &x_2 &\cdots &x_{n+3} \\
-x_1    & 0      &y_{23} &\cdots & y_{2,n+3}  \\
-x_2    &-y_{23} &0 &\cdots & y_{3,n+3}\\
\vdots  &0       &\vdots  &0 &  \vdots \\
-x_{n+3} &-y_{2,n+3} &-y_{3,n+3} &\cdots  & 0
\end{array}\! \right]
\]
To see why this last statement holds, consider any $5 \times 5$ 
submatrix of the form
\[
\left[ \!
\begin{array}{cccccc}
0        &x_{i_1} &x_{i_2} &x_{i_3}  &x_{i_4} \\
-x_{i_1} & 0      &*       &*        & *  \\
-x_{i_2} & *      &0       &*        & *\\
-x_{i_3} & *      &*       &0        &  * \\
-x_{i_4} & *      &*       &*        & 0
\end{array}\! \right]
\]
A computation shows that if $P$ is the $4 \times 4$ Pfaffian of
the block $*$, then $x_{i_1}\cdot P \in I$; since $I$ is
prime and nondegenerate this implies $P \in I_2$.
\vskip .1in
\noindent For part $(3)$, it follows from $(2)$ that 
if $I$ has a linear $n^{th}$ syzygy of rank $n+3$, then $I_2$ contains $\mbox{Pfaff}_4(N)$ of
a skew-symmetric matrix $N$, so we have polynomials 
of the form
\[
y_{ij}x_k - y_{ik}x_j + y_{jk}x_i \in I_2. 
\]
Since $I$ is a semigroup ideal, every term in such a 
polynomial has the same weight with respect to the
semigroup. As a semigroup ideal, $I$ is generated by 
binomials, so one of the $y_{ij}$ must vanish, because otherwise 
subtracting $y_{ij}x_k - y_{ik}x_j \in I_2$ from 
$y_{ij}x_k - y_{ik}x_j + y_{jk}x_i$ shows that 
$y_{jk}x_i \in I_2$, contradicting nondegeneracy of $I$.
Repeating this shows that in the semigroup case the matrix $N$ is of the form
\[
\left[ \!
\begin{array}{cc}
0       &M \\
-M^t    & *
\end{array}\! \right],
\]
where $M$ is a $1$--generic $v \times w$ matrix and $v, w \ge 2$. It is
easy to show that 
$I_2(M) \subseteq \mbox{Pfaff}_4(N)$, and the result follows.
\section{Toric specializations of Rees algebras}
\subsection{Rees algebras of Koszul cycles}
The skew-symmetric matrices which arose in \S 2 in conjunction
with the linear first syzygies are most naturally  studied in the
setting of Rees algebras. Let $P = \k[x_1,\ldots,x_d]$, and let 
$K_\bullet$ denote the Koszul complex on $\{x_1,\ldots,x_d\}$
\[
K_\bullet : \mbox{ } 0 \rightarrow \bigwedge^d(\k^d)\otimes P(-d)
\stackrel{\delta_d}{\rightarrow} \bigwedge^{d-1}(\k^d)\otimes P(-d+1)
\stackrel{\delta_{d-1}}{\rightarrow} \bigwedge^{d-2}(\k^d)\otimes P(-d+2)
\stackrel{\delta_{d-2}}{\rightarrow} \cdots
\]
Let $C_i = \im(\delta_i)$ be the module of $i^{th}$ cycles in $K_\bullet$,
and put
\[
\mathcal{I}_i = \{I \subseteq \{1,\ldots,d\} \mid |I|=i \}.
\]
Then the symmetric algebra on the free module $K_i$ is:
\[
S(K_i) = P[y_I \mid I \in \mathcal{I}_i]
\]
The presentation 
\[
K_{i+1} \longrightarrow K_{i} \longrightarrow C_{i}\longrightarrow 0
\]
gives a presentation for the symmetric algebra of $C_i$. Let $J_i = \langle z_I \mid I \in \mathcal{I}_{i+1} \rangle$,
where if $I = \{ a_1 < a_2 < \cdots < a_{i+1} \}$ then 
\[
z_I = \sum \limits_{j=1}^{i+1} (-1)^{j+1} x_{a_j}y_{I \setminus a_j} \mbox{ and } S(C_i) \simeq S(K_i)/J_i.
\]
For $i=d-2$, any $I \in \mathcal{I}_{d-1}$ is of the form 
$[d] \setminus k$, so $I \setminus a_j$ has the form
\[
I \setminus a_j = [d] \setminus \{k, j \} = \widehat{jk} \mbox{ for some }j,k.
\]
Thus, for $i=d-2$, the $z_I$ are exactly the elements denoted $\partial(\omega)$
in Equation~\ref{schreyerEq}, and the reason for the choice of notation in \S 2.
The Rees algebra $R(C_i)$ is $S(C_i)$ modulo the $P$-torsion.
As noted earlier, for $i=d-2$ these algebras were first
investigated in \cite{hu}, and the free resolution of $R(C_{d-2})$
and $S(C_{d-2})$ over $S(K_{d-2})$ was determined by Kustin in \cite{k86}. In \cite{htz},
Herzog, Tang and Zarzuela studied properties of Gr\"obner and
sagbi bases for $R(C_i)$ and $S(C_i)$, concentrating on the cases $i=2$ and
$d-2$; they also conjectured that $R(C_i)$ is Cohen-Macaulay for 
all $i$, and that the $P$-torsion of $S(C_i)$ could be described
simply as $0:_{S(C_i)}x_j$ for any $j \in \{1, \ldots,d\}$.

In \cite{w}, Weyman used the geometric method of computing
syzygies \cite{w2} to prove these conjectures. In fact, Weyman shows
that $R(C_i)$ is a normal, Cohen-Macaulay domain, with
rational singularities, and gives a representation theoretic 
description of the free resolution of $R(C_i)$ over $S(K_i)$. 
In different language, \cite{vB2} studies $Proj(S(C_i))$, 
obtaining results on rank two bundles on curves.

As noted in \S 2, the $z_I$ are
quadrics of high rank, so the question is how to specialize
the $z_I$ so that they have rank at most four, but where
the specialized ideal remains prime. 

\subsection{Toric specializations of $J_i$}
Let $\Delta$ be an $n$--dimensional simplicial complex on $d$ vertices.
We associate to $\Delta$ an ideal $J_\Delta$ which is a specialization of
$J_{d-n-1}$. If $\Delta$ is a pseudomanifold, then $J_\Delta$ is toric; as
a pseudomanifold $\Delta$ has a natural top homology class, which 
corresponds to a minimal syzygy on $J_\Delta$. The specializations 
in \S 2 coming from cycles are instances of this construction;
motivated by this, we focus on specializations for which
the underlying simplicial complex is a pseudomanifold.

\begin{defn}
An oriented pseudomanifold $\Delta$ of dimension $n$ on $d$ vertices consists of a set
of oriented $(n+1)$-subsets of $\{1,\ldots,d\}$ such that
\begin{enumerate}
\item{each $n$-subset of $\{1,\ldots,d\}$ is contained in 
exactly zero or two elements of $\Delta$; in the latter case 
the orientations must cancel.}

\item{$\Delta$ is strongly connected: the dual graph $G(\Delta)$ is
 connected, where $G(\Delta)$ has a vertex for each $n$--face of
 $\Delta$, and two vertices are joined by an edge if the corresponding
$n$--faces share an $(n-1)$--face.}
\end{enumerate}
\end{defn}

\begin{defn}\label{PMideal}
Let $\Delta$ be an $n$--dimensional oriented pseudomanifold on 
$d$ vertices. $J(\Delta)$ is the specialization
of $J_{d-n-1}$ obtained by setting $y_{I \setminus a_j} =0$ if 
$ \widehat{I \setminus a_j} \not \in \Delta$. So $y_{I \setminus a_j}=0$ iff $I \setminus a_j$ 
is in the Alexander dual $\Delta^*=\{ \widehat{\gamma}  \mid \gamma \not\in \Delta \}$.
\end{defn}

\begin{exm}\label{exm:bdrySimplex}
If $\Delta$ is the boundary of an $k$--simplex, then
$\Delta$ has $k+1$ vertices and $k+1$ faces
of dimension $k-1$, and $J_\Delta$ is the ideal of 2 by 2 
minors of a generic 2 by (k+1) matrix. For $k=2$, we have
$d=3$, $n=1$. Then $J_1$ is the ideal
\[
\begin{array}{ccccc}
z_{12} & = & x_1y_2-x_2y_1 & = x_1y_{\widehat{13}}-x_2y_{\widehat{23}}\\
z_{13} & = & x_1y_3-x_3y_1 & = x_1y_{\widehat{12}}-x_3y_{\widehat{23}}\\
z_{23} & = & x_2y_3-x_3y_2 & = x_2y_{\widehat{12}}-x_3y_{\widehat{13}}
\end{array}
\]
yielding the $2 \times 2$ minors of the $2 \times 3$ matrix appearing in 
Example~\ref{exm:3vars}. We relabel $y_I$ as $y_{\widehat{[d]\setminus I}}$ to make the 
connection to the cycle more intuitive. 
\end{exm}
\begin{lem}\label{specialLem}
With the relabelling introduced above, 
\[
J(\Delta) \simeq \langle x_i y_{\widehat{\sigma}} - x_j y_{\widehat{\tau}} \mid \sigma, \tau \in \Delta_n \mbox{ satisfy } \sigma \setminus \{i\} = \tau \setminus \{j\} \rangle.
\]
\end{lem}
\begin{proof}
The simplicial coboundary map $\partial_n^*: C^{n-1}(\Delta) \rightarrow C^{n}(\Delta)$ has two
nonzero entries in each column, $+1$ and $-1$. Index the source and
target of $\partial_n^*$ by complementary faces, so that $C^{n-1}(\Delta)$ is indexed by $\mathcal{I}_{d-n}$ and $C^{n}(\Delta)$ by $\mathcal{I}_{d-n-1}$. 
Because Koszul and simplical cohomology agree, choosing an oriented 
basis for the Koszul cycles to agree with the orientation of 
$\Delta$ and arbitrary orientations for cycles not in $\Delta$ yields the result. 
\end{proof}
In concrete terms, the ideal $J_{d-n-1}$ consists of column sums of a ${d \choose d-n-1} \times {d \choose d-n}$ 
matrix, with entry $(J,I) \in \mathcal{I}_{d-n-1}\times\mathcal{I}_{d-n}$ zero if $J \not\subseteq I$, and entry
$(-1)^{j+1}x_{a_j}y_{I\setminus a_j}$ if $I = [a_1,\ldots,a_{d-n}]$ and $J = [a_1,\ldots,\widehat{a_j},\dots,a_{d-n}]$.
Lemma~\ref{specialLem} replaces $(-1)^{j+1}$ with the sign of $\widehat{I}$ in the boundary of $\widehat{J}$,
and $y_{J}$ with $y_{\widehat{[d]\setminus J}}$.
\begin{exm}\label{exm:6cycleagain}
We revisit Example~\ref{exm:6vars}, where $\Delta$ is the six-cycle with
orientation $\{[i,i+1], i \in \{1,\ldots,5\}, [6,1]\}$.
 Since $d=6$ and $n=1$, we consider $J_{d-n-1}=J_4$.
Taking lex ordered bases for the Koszul classes yields, for example:
\[
\begin{array}{ccc}
z_{23456} & = & x_2y_{3456}-x_3y_{2456}+x_4y_{2356}-x_5y_{2346}+x_6y_{2345}\\
         & = & x_2y_{\widehat{12}}-x_3y_{\widehat{13}}+x_4y_{\widehat{14}}-x_5y_{\widehat{15}}+x_6y_{\widehat{16}}\\
         & \leadsto & x_2y_{\widehat{12}}+x_6y_{\widehat{16}}.\\
\end{array}
\]
Modifying the basis for Koszul classes as in Lemma~\ref{specialLem} changes the last
expression to $x_2y_{\widehat{12}}-x_6y_{\widehat{16}}.$ This explains the
linear change of variables $y_{1d}\mapsto -y_{1d}$ in the proof of Theorem 1.6.
\end{exm}
$J(\Delta)$ is typically not prime, but if 
$J_\Delta = J(\Delta):\langle \prod_{i=1}^dx_i \rangle ^{\infty}$, then 
Lemma~\ref{lem:JDprime} shows that $J_\Delta$ is prime, toric, 
and is the defining ideal of a Rees algebra.

\begin{exm}\label{exm:S2}
Consider the two-dimensional pseudomanifold which is a triangulation
of $S^2$, obtained by coning over the boundary of the complex below 
with a seventh vertex.
\begin{figure}[h]
\epsfig{file=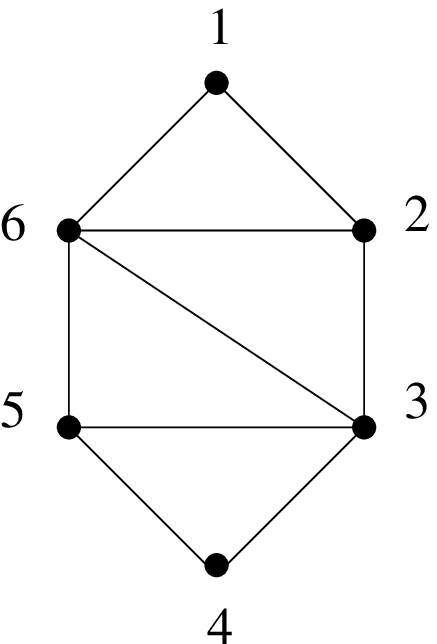,height=1.5in,width=1.5in}
\end{figure}

\noindent The resulting ideal $J(\Delta)$ is generated by fifteen
quadrics; $J_\Delta$ defines a toric seven--fold in $\p^{16}$, which is
Cohen-Macaulay of degree $73$, and has graded betti numbers:
\begin{small}
$$
\vbox{\offinterlineskip 
\halign{\strut\hfil# \ \vrule\quad&# \ &# \ &# \ &# \ &# \ &# \
&# \ &# \ &# \ &# \ &# \ &# \ &# \
\cr
total&1&21&163&447&631&575&377&168&47&8\cr
\noalign {\hrule}
0&1 &--&--& --&--& --& --&--& --& --  \cr
1&--&17 &19&1&--&--&   --&--& --& --    \cr
2&--&4 &144 & 444&500&209&8&--& --&--     \cr
3&--&-- &-- &2&131&365& 333&93& 2&  --    \cr
4&--&-- &-- &--&--&1& 36&75& 45& 8     \cr
\noalign{\bigskip}
\noalign{\smallskip}
}}
$$
\end{small}
\end{exm}
\noindent The existence of a linear second syzygy (of rank
seven) follows from Theorem~\ref{thm:third}, which we prove next.
\begin{defn}\label{toricIdeal}Let 
\[
\mathcal{A} = \{ {\bf a}_1,\ldots, {\bf a}_q \} \subseteq \mathbb{Z}^p
\]
and let $A$ be the matrix with $i^{th}$ column ${\bf a}_i$. The
toric ideal 
\[
I_{\mathcal{A}} = \langle x^{\bf \alpha}- x^{\bf \beta} \mid \alpha, \beta \in  \mathbb{N}^q \mbox{ and } \alpha-\beta \in ker(A)\rangle.
\]
To any ${\bf m} \in \mathbb{Z}^p$ we associate a simplicial complex
\[
\Delta_{{\bf m}}(\mathcal{A})=\{I \subseteq \{1,\ldots,q\} \mid
{\bf m}-\sum\limits_{i\in I}{\bf a}_i \mbox{ lies in }\mathcal{A} \}.
\]
\end{defn}
\noindent A result of Hochster shows that 
\[
\widetilde{H}_j(\Delta_{{\bf m}}(\mathcal{A}),\k) = Tor_j(I_\mathcal{A},\k)_{{\bf m}}.
\]
For a proof and more details on $I_\mathcal{A}$ see Chapter 9 of \cite{ms}.
Let $\Delta$ be an oriented $n$--dimensional pseudomanifold on $d$--vertices,
and for each $i \in \{1,\ldots,d \}$, 
assign the vertex $v_i$ the weight ${\bf e}_i \in \mathbb{Z}^{d+1}$, 
and for a facet $\sigma = \{ v_{i_1}, \ldots, v_{i_{n+1}} \} \in \Delta_n$, 
\[
wt(\sigma)={\bf e}_0-\sum\limits_{j=1}^{n+1} {\bf e}_{i_j}.
\]
For the remainder of this section, 
\[
\mathcal{A} = \{{\bf e}_1, \ldots, {\bf e}_d, wt(\sigma) \mid \sigma \in \Delta_n \} 
\subseteq \mathbb{Z}^{d+1}.
\]
\begin{lem}\label{lem:homotopic}
Let $\Delta$ be an oriented $n$--dimensional 
pseudomanifold on $d$--vertices, with top homology class
$H_n(\Delta)$. Then $\Delta$ is homotopic to $\Delta_{{\bf e}_0}(\mathcal{A})$
\end{lem}
\begin{proof}
Consider the complex $\Delta_{{\bf e}_0}(\mathcal{A})$. Let $\tau = \{i_1,\ldots,i_k \}$
define a $(k-1)$--face of $\Delta$. Since $\Delta$ is a pseudomanifold,
there is an $n$--simplex $\{i_1,\ldots,i_k, i_{k+1},\ldots, i_{n+1}
\}$ containing $\tau$, so
\begin{equation}\label{eqn:1eq}
{\bf e}_0-{\bf e}_{i_1}-\cdots -{\bf e}_{i_k} = 
{\bf e}_0-\sum\limits_{j=1}^{n+1}{\bf e}_{i_j} +
\sum\limits_{j=k+1}^{n+1}{\bf e}_{i_j},
\end{equation}
hence $\tau  \in \Delta_{{\bf e}_0}(\mathcal{A})$, which
implies that $\Delta \subseteq  \Delta_{{\bf e}_0}(\mathcal{A})$. 
On the other hand, 
\begin{equation}\label{eqn:2eq}
{\bf e}_0-({\bf e}_0-{\bf e}_{i_1}-\cdots -{\bf e}_{i_k}) = 
\sum\limits_{j=1}^{k}{\bf e}_{i_j} \in \mathcal{A},
\end{equation}
so the cone over every $(k-1)$--face of $\Delta$ is also in 
$\Delta_{{\bf e}_0}(\mathcal{A})$. This is most easily visualized as
adding a cone vertex (corresponding to the variable $y_\sigma$) 
over every $n$--face $\sigma$ of $\Delta$. Any $\sum_{i \in I}{\bf a}_i$ such 
that $I \in \Delta_{{\bf e}_0}(\mathcal{A})$ can have at most one ${\bf a}_i$ with
${\bf e}_0$ coefficient one. If there are no such ${\bf a}_i$, then $I$
must correspond to an element of $\Delta_n$. Otherwise, let $k \in I$ with
${\bf a}_k = {\bf e}_0 - \sum_{j \in \sigma} {\bf e}_j$. Then
\[
{\bf e}_0 - \sum_{i \in I} {\bf a}_i = {\bf e}_0 - ({\bf e}_0-\sum\limits_{j \in \sigma}{\bf e}_j) - (\sum\limits_{i \in I' = I \setminus k}{\bf a}_i) \in \mathcal{A},
\]
which can only occur if $I' \subseteq \sigma$. 
Thus, all faces of $\Delta_{{\bf e}_0}(\mathcal{A})$
are described by Equations \ref{eqn:1eq} and \ref{eqn:2eq}. 
Since each cone over a face $\sigma$ (the faces appearing in
Equation~\ref{eqn:2eq}) is homotopic to the face $\sigma$, the result follows.
\end{proof}

\begin{defn}
Let $L_\Delta \subseteq S=\k[x_1,\ldots,x_d,\{ y_\sigma  \mid \sigma \in \Delta_n \}]$ be the defining ideal of 
the Rees algebra $P[I_{\Delta^*} t]$, where $I_{\Delta^*}$ is the Stanley-Reisner ideal of the Alexander dual of $\Delta$, and the map $S \rightarrow P[I_{\Delta^*} t]$ is defined via $y_\sigma \mapsto t\prod_{i \not \in \sigma}x_i$.
\end{defn}

\begin{lem}\label{lem:JDprime}
Let $\Delta$ be an oriented $n$--dimensional 
pseudomanifold on $d$--vertices. Then $J_\Delta = L_\Delta = I_{\mathcal A}$ is a prime toric ideal.
\end{lem}
\begin{proof}
Since the Stanley-Reisner ideal $I_{\Delta^*}$ of the Alexander dual of $\Delta$ is
\[
I_{\Delta^*} = \langle \prod\limits_{i \not\in \sigma} x_i  \mid \sigma \in max(\Delta)\rangle,
\]
the following inclusions are immediate:
$$J(\Delta) \subseteq L_\Delta \subseteq I_{\mathcal A}.$$
Let $f=\prod\limits_{i=1}^d x_i$. Since $\Delta$ is strongly connected, we have
$$J(\Delta) S_f = I_{\mathcal A} S_f,$$
Therefore
$$J(\Delta) S_f = L_\Delta S_f =  I_{\mathcal A} S_f.$$
By contracting back to $S$,
noting that $L_\Delta$ and $I_{\mathcal A}$ are prime, and $J(\Delta) S_f
\cap S = J_\Delta$, we obtain the desired equalities.
\end{proof}

\noindent Combining Lemma~\ref{lem:JDprime} with Lemma~\ref{lem:homotopic} yields
Theorem~\ref{thm:third}. In Example~\ref{exm:S2}, the linear syzygy 
constructed from Theorem~\ref{thm:third} generates $Tor_3(S/J_\Delta,\k)_{4}$, 
but in general its dimension can be larger.
\begin{prop}\label{prop:Biglins}
Let $\Delta$ be an oriented $n$--dimensional 
pseudomanifold on $d$--vertices, with top homology class
$H_n(\Delta)$. Suppose $\Delta' \subseteq \Delta$ is a bipyramid on
vertices $x_1,x_2$ over a $n\!-\!1$--dimensional complex $\Delta''$.
Then if $|\Delta''_{n-1}| = k$, $J_\Delta$ contains the 
two by two minors of a matrix of the form:
\[
\left[ \!
\begin{array}{ccccc}
x_1 & y_1 &y_2 & \cdots &y_k\\
x_2 & y_1'&y_2'& \cdots & y_k'
\end{array}\! \right]
\]
\end{prop}
\begin{proof}
Let  $\Delta''_{n-1} = \{\sigma_1,\ldots,\sigma_k \}$. Each
$\sigma_i$ yields a pair of $n$--faces of $\Delta$: 
the cone of $\sigma$ with $x_1$ and $x_2$. Associate to the cone 
of $\sigma_i$ with $x_2$ the variable $y_i$, and to the cone 
of $\sigma_i$ with $x_1$ the variable $y_i'$. Then 
\[
\deg(x_1y_i') = {\bf e}_1 + ({\bf e}_0-\sum\limits_{i \in \sigma}{\bf e}_i -{\bf e}_1) =
{\bf e}_2 + ({\bf e}_0-\sum\limits_{i \in \sigma}{\bf e}_i -{\bf e}_2) =\deg(x_2y_i).
\]
A similar argument shows that $\deg(y_iy_j')=\deg(y_jy_i')$.
\end{proof}
\noindent The point is that if $\dim_\k Tor_{n}(S/J_\Delta,\k)_{n+1}=1$, then 
$\Delta$ cannot contain a bipyramid over an $n-2$--dimensional
$\Delta''$ with $|\Delta''_{n-2}| \ge n$. By Example~\ref{exm:bdrySimplex},
this also implies that $\Delta$ is not the boundary of an 
$n$--simplex.
\begin{que} 
What are necessary and sufficient conditions on
a pseudomanifold $\Delta$ so that $\dim_\k Tor_n(S/J_\Delta,\k)_{n+1}=1$? What
are necessary and sufficient conditions on a pseudomanifold $\Delta$ so that
Conjecture~\ref{conj:GG} holds for $S/J_\Delta$?
\end{que}

\noindent Any one--dimensional $\Delta$ which is a cycle on five or more 
vertices gives $J_\Delta$ with $\dim_\k Tor_2(S/J_\Delta,\k)_{3}=1$. 
For the toric varieties produced from two--dimensional 
pseudomanifolds on at most six vertices, there are no counterexamples 
to Conjecture~\ref{conj:GG}.
However, Example~\ref{exm:S2} may be specialized to yield 
a toric counterexample:
\noindent \begin{exm}
Identifying the vertices labeled $1$ and $4$ in Example \ref{exm:S2}
and saturating yields an ideal 
which defines a toric six--fold in $\p^{15}$, which is
Cohen-Macaulay of degree $56$, and has graded betti numbers:
\begin{small}
$$
\vbox{\offinterlineskip 
\halign{\strut\hfil# \ \vrule\quad&# \ &# \ &# \ &# \ &# \ &# \
&# \ &# \ &# \ &# \ &# \ &# \ &# \
\cr
total&1&25&177&549&816&676&449&255&67&5\cr
\noalign {\hrule}
0&1 &--&--& --&--& --& --&--& --& --  \cr
1&--&19 &30&1&--&--&   --&--& --& --    \cr
2&--&6 &147 & 546&788&484&45&--& --&--     \cr
3&--&-- &-- &2&28&192&404&255& 64& 3    \cr
4&--&-- &-- &--&--&--& --&--& 3& 2     \cr
\noalign{\bigskip}
\noalign{\smallskip}
}}
$$
\end{small}
This specialization does not correspond to a pseudomanifold;
there is an edge of $\Delta$ (connecting vertex $1=4$ and vertex 7),
which lies on four two-faces.
\end{exm}
\begin{que} If $\Delta$ is a pseudomanifold, which specializations of
$J_\Delta$ are toric? If $\dim_\k Tor_n(S/J_\Delta,\k)_{n+1}=1$, which
specializations of $J_\Delta$ are toric, {\em and} preserve uniqueness
of the top linear syzygy?
\end{que}
\section{Examples with $2LP=2$}
In this final section, we consider geometric reasons for the failure
of Conjecture~\ref{conj:GG}. We restrict our attention to the case 
where $2LP=2$, and focus on two concrete classes of examples: 
curves and toric surfaces. In Example~\ref{exm:5vars}, the 
divisor which was used to embed
the curve was special. The next result shows that this is not
an isolated phenomenon. Recall that normal generation means $S/I_C$ is projectively normal,
and normal presentation means $I_C$ is generated by quadrics.
\begin{prop}\label{prop:curve1}
If $D$ is a very ample divisor on $C$ such that $C \subseteq \p(H^0(D)^*)$ is
normally generated and normally presented, and $b_{23}(S/I_C)=1$, then 
$D$ is special.
\end{prop}
\begin{proof}
Let $r+1 = h^0({\mathcal O}_C(D))$ and $d = deg(D)$. 
The assumption on normal generation means that  $h^1({\mathcal I}_C(t))=0$
 for all $t\ge 0$, yielding the exact sequence
\vskip .05in
\[
0 \longrightarrow H^0({\mathcal I}_C(t))\longrightarrow H^0({\mathcal O}_{\p^r}(t))
\longrightarrow H^0({\mathcal O}_{C}(t))\longrightarrow 0.
\]
\vskip .1in
\noindent If $D$ is nonspecial, then $h^1(({\mathcal O}_{C}(tD))=0$ for all $t
\ge 1$, which for $t=1$ implies:
\[
g=d-r.
\]
For $t=2$, Riemann-Roch and the short exact sequence above 
shows that
\[
h^0({\mathcal I}_C(2))={r+2 \choose 2}-2d-1+g.
\]
Since $I_C$ is generated by quadrics, having a single linear syzygy implies that
\[
h^0({\mathcal I}_C(3)) = (r+1) \cdot h^0({\mathcal I}_C(2)) - 1 = {r+3 \choose 3}-3d-1+g.
\]
Eliminating $g$ from these equations shows that
\[
d = \frac{r^3-r-3}{3(r-1)},
\]
which has no integral solutions. 
\end{proof}
\noindent Taking this as our cue, we study nonspecial $D$, 
such that $C \subseteq \p(H^0(D)^*)$ is 
normally generated and presented, with $2LP(S/I_C)=2$.
When $\deg(D) \ge 2g+1+a$, a result of Schreyer 
(\cite{e2}, Theorem 8.17) shows that $2LP(S/I_C) \ge a+\lfloor 
\frac{g}{2} \rfloor$, which is forced due to a factorization of
$D$. This implies Conjecture \ref{conj:GG} holds if 
$2LP(S/I_C)=2$ and $\deg(D) \ge 2g+1+a$. Slightly modifying 
Schreyer's proof yields:
\begin{lem}\label{lem:curve2}
If $\mathcal{L}$ is a very ample line bundle on $C$ such that 
\[
\deg \mathcal{L} = g + \Big\lceil \frac{g}{2} \Big\rceil +a, \mbox{ } a \ge 1,
\]
then $2LP(S/I_C) \ge a-1$.
\end{lem}
\begin{proof}
Brill--Noether implies that if $m = 1+ \lceil \frac{g}{2} \rceil$,
then $C$ carries a $g^1_m$. Let $D$ be a divisor in this system. 
Then
\[
h^0(\mathcal{L}(-D)) = g+ \Big\lceil \frac{g}{2} \Big\rceil +a - (1+ \Big\lceil
\frac{g}{2} \Big\rceil) +1 -g + h^1(\mathcal{L}(-D)).
\]
If $X$ is a curve embedded by a complete linear series $\mid \mathcal{A}
\mid$ and $B$ is a divisor on $X$ such that $h^0(B)\ge s+1$, then 
if $h^0(\mathcal{A}(-B))\ge t+1$, $I_X$ contains the $2 \times 2$ minors 
of a $(s+1) \times (t+1)$ matrix. Theorem 8.12 of \cite{e2} implies that
\[
2LP(S/I_X) \ge s+t-1.
\]
Applying this to $\mathcal{L}$ and $D$ yields the lemma. If  
$\mathcal{L}(-D)$ is special, then in fact $2LP(S/I_C) \ge a$.
\end{proof}
\begin{lem}\label{lem:curve3}
If $D$ is a very ample, nonspecial divisor on $C$ of degree $d$, 
such that $C \subseteq \p(H^0(D)^*)$ is
normally generated and normally presented and
$2LP(S/I_C)=2$, then the graded betti numbers of $S/I_C$ are given by:
\begin{small}
$$
\vbox{\offinterlineskip 
\halign{\strut\hfil# \ \vrule\quad&# \ &# \ &# \ &# \ &# \ &# \
&# \ &# \ &# \ &# \ &# \ &# \ &# \
\cr
total & 1& $b_2$ & $\cdots$ & $\cdots$  & $\cdots$ & $b_{r+1}$ &\cr
\noalign {\hrule}
0&1 &--  &--  &--  &--     & --      &\cr
1&--&$b_2$ &$b_3$ &--  &--     & --      &\cr
2&--&--  &$b_4$ &$b_5$ &$\cdots$ & $b_{r+1}$ &\cr   
\noalign{\bigskip}
\noalign{\smallskip}
}}
$$
\end{small}
where 
\[ 
\begin{array}{ccc}
b_2 & = & {d-g \choose 2}-g\\
b_3 & = & (d-g-1)({d-g-1 \choose 2}-g) -{d-g-1\choose 3}\\
\end{array}
\]
and when $i \ge 4$,
\[
b_i  =  {d-g-1 \choose i}-(d-g-1) {d-g-1 \choose i-1} +
g{d-g-1 \choose i-2}
\]
\end{lem}
\begin{proof}
The hypothesis of projective normality means that all values of the Hilbert function
of $S/I_C$ can be computed from $H^0(D)$. As in the proof of 
Proposition~\ref{prop:curve1}, Riemann-Roch and the hypothesis that $D$ is
nonspecial yield these values. The assumption that $I_C$ is
generated by quadrics and that $2LP(S/I_C)=2$ means that there are no
overlaps in the resolution, hence the Hilbert function in fact
determines the resolution.
\end{proof}

Since the $b_i$ are positive, the hypotheses of Lemma~\ref{lem:curve3} impose very strong constraints
on $D$: for curves of genus at most six, the only possible values for genus and degree are:
\begin{center}
\begin{supertabular}{|c|c|c|c|c|c|c|c|}
\hline genus &  $0$ & $1$ & $2$ & $3$ & $4$ & $5$ & $6$ \\
\hline degree&  $3$ & $5$ & $6$ & $8$ & $9$ & $11$ & $12$ \\
\hline
\end{supertabular}
\end{center}
\vskip .05in
Lemma~\ref{lem:curve2} implies that Conjecture~\ref{conj:GG} holds for
all such pairs. For a curve of genus seven, a divisor satisfying the 
hypotheses of Lemma~\ref{lem:curve3} must have degree $13$ or $14$. 
In the latter case, Conjecture~\ref{conj:GG} again holds by 
Lemma~\ref{lem:curve2}.
\noindent \begin{exm}
On a curve of genus seven, a divisor $D$ of degree $13$ satisfying the 
hypotheses of Lemma~\ref{lem:curve3} has graded betti numbers
\vskip .01in
\begin{small}
$$
\vbox{\offinterlineskip 
\halign{\strut\hfil# \ \vrule\quad&# \ &# \ &# \ &# \ &# \ &# \
&# \ &# \ &# \ &# \ &# \ &# \ &# \
\cr
total & 1& 8 & 30 & 46  & 30 & 7 &\cr
\noalign {\hrule}
0&1 &--  &-- &--  &--  & --      &\cr
1&--&8   &5  &--  &--  & --      &\cr
2&--&--  &25 & 46 & 30 & 7       &\cr   
}}
$$
\end{small}
\noindent For a general curve, $D = K+p_1+p_2+p_3+p_4-q_1-q_2-q_3$ 
has such betti numbers. A computation shows that all linear syzygies have 
rank $\ge 5$, and  that $I_C$ can be generated by quadrics of rank $\le 4$. 
So Conjecture~\ref{conj:GG} can fail even if $D$ satisfies the hypotheses of
Lemma~\ref{lem:curve3}.
\end{exm}

Nevertheless, there are classes of objects where the constraints of
Lemma~\ref{lem:curve3} are strong enough to prove Conjecture~\ref{conj:GG}. 
Let $X$ be a complete toric surface and $D$ a very ample divisor on $X$.
Then $H^0(D)$ corresponds to the set of integral points of a lattice 
polygon $P$, and $S/I_X$ is Cohen-Macaulay and three-regular, so 
we may repeat the analysis above. Let $v(P)$ denote the volume of $P$, 
$\partial(P)$ the number of lattice points on the boundary 
of $P$, and $i(P)$ the number of interior lattice points 
of $P$; write $X_P$ for the projective embedding 
of $X$ determined by $P$.
\begin{prop}\label{prop:surface}
If $X_P$ is a toric surface generated by quadrics and 
$2LP(S/I_{X_P})=2$, then Conjecture \ref{conj:GG} holds.
\end{prop}
\begin{proof}
Pick's theorem shows that 
for a divisor $D$ on a toric surface, 
\[
D^2 = 2v(P) = 2i(P)+\partial(P)-2.
\]
The surface $X_P$ is projectively normal \cite{BGT}, so 
if $C = X_P \cap H$ for a general hyperplane section, the graded betti 
numbers of $C$ and $X$ are identical. Riemann-Roch and 
adjunction show that the genus of $C$ is equal to $i(P)$, and so $C$ is embedded by
a divisor of degree $2g + 1 + \partial(P)-3.$
Slicing with a hyperplane as in \cite{sch} and applying 
Green's theorem \cite{g2} shows that 
\[
Tor_i(S/I_C,\k)_{i+2}=0 
\]
for all $i \le \partial(P)-3$; in \cite{k}, Koelman shows that 
if $Tor_1(S/I_C,\k)_{3}=0$ and $Tor_2(S/I_C,\k)_{4}\ne 0$, then 
$\partial(P)=4$. By Lemma \ref{lem:curve2} if $\mathcal{L}$ is
a line bundle on $C$ of degree greater than $g +\lceil \frac{g}{2}
\rceil+4$, then $2LP(S/I_C)\ge 3$, so the assumption that
$2LP(S/I_C)=2$ implies that 
\[
2g + 1 + \partial(P)-3 = \deg(D|_D) \le g +\Big\lceil \frac{g}{2}
\Big\rceil+3,
\]
hence $\big\lfloor \frac{i(P)}{2}\big\rfloor + \partial(P) \le 5$. This 
means either $\partial(P)=4$ and $i(P) \in \{0,\ldots,3\}$ or 
$\partial(P)=5$ and $i(P)$ is $0$ or $1$. There are only finitely
many such polygons, and a check verifies the conjecture. 
\end{proof}

\noindent {\bf Acknowledgments and Remarks}:  Macaulay2 computations were
essential to our work. We also thank MSRI and BIRS, where portions 
of this work took place, and the referee for useful comments. 

\bibliographystyle{amsalpha}

\begin{thebibliography}{10}
\bibitem{BGT} W. ~Bruns, J. ~Gubeladze, N.V. ~Trung,
        {\em Normal polytopes, triangulations, and Koszul algebras},
         J. Reine Angew. Math.  {\bf 485}  (1997), 123--160.

\bibitem{e1} D. ~Eisenbud,
        {\em Green's conjecture: an orientation for algebraists},
        in ``Free resolutions in commutative algebra and algebraic geometry''
        (Sundance, UT, 1990),  51--78, Res. Notes Math., 2, Jones and
        Bartlett, Boston, MA, 1992

\bibitem{e2} D. ~Eisenbud,
        {\em The geometry of syzygies},
        Graduate Texts in Mathematics, vol.~229,
        Springer-Verlag, Berlin-Heidelberg-New York, 2005.

\bibitem{e3} D. ~Eisenbud,
        {\em Linear sections of determinantal varieties},
Amer. J. Math.  {\bf 110}  (1988),  no. 3, 541--575. 

\bibitem{ek} D. ~Eisenbud, J.~Koh,
{\em Some linear syzygy conjectures},
Adv. Math. {\bf 90} (1991), 47--76.

\bibitem{eks} D. ~Eisenbud, J.~Koh, M. ~Stillman, 
{\em Determinantal equations for curves of high degree},
Amer. J. Math. {\bf 110} (1988), 513--539.

\bibitem{g1} M. ~Green,
{\em Quadrics of rank four in the ideal of a canonical curve},
Invent. Math. {\bf 75} (1984), 85--104. 

\bibitem{g2} M. ~Green,
{\em Koszul cohomology and the geometry of projective varieties},
J. Differential Geom. {\bf 19} (1984), 125--171.

\bibitem{htz} J. ~Herzog, Z. ~Tang, S. ~Zarzuela, 
{\em Symmetric and Rees algebras of Koszul cycles and their Gr\"obner
  bases},
Manuscripta Math. {\bf 112} (2003), 489--509.
 
\bibitem{hu} C. ~Huneke, B. ~Ulrich, 
{\em Divisor class groups and deformations},
Amer. J. Math. {\bf 107} (1986), 1265--1303.

\bibitem{k} R.~Koelman, {\em A criterion for the ideal of a projectively
    embedded toric surface to be generated by quadrics}, Beitr\"age
  Algebra Geom. \textbf{34} (1993), 57--62.

\bibitem{ks} J.~Koh, M. ~Stillman, 
{\em Linear syzygies and line bundles on an algebraic curve},
J. Algebra  {\bf 125}  (1989), 120--132. 

\bibitem{k86} A. ~Kustin,
{\em The minimal free resolutions of the Huneke-Ulrich deviation two Gorenstein ideals},
J. Algebra  {\bf 100} (1986), 265--304. 

\bibitem{ms} E. ~Miller, B. ~Sturmfels,
        {\em Combinatorial commutative algebra},
        Graduate Texts in Mathematics, vol.~227,
        Springer-Verlag, Berlin-Heidelberg-New York, 2005.

\bibitem{sch} H.~Schenck, {\em Lattice polygons and Green's theorem}, 
Proc.\  Amer.\ Math.\ Soc.\ \textbf{132} (2004) 3509--3512.

\bibitem{s} F-O. ~Schreyer,
              {\em A standard basis approach to syzygies of canonical
              curves},
              J. Reine Angew. Math.  {\bf 421}  (1991), 83--123.
\bibitem{vB2} H. ~von Bothmer,
{\em Generic syzygy schemes},
J. Pure App. Alg., {\bf 208} (2007), 867--876.

\bibitem{w} J. ~Weyman, 
{\em Application of the geometric technique of calculating syzygies to
  Rees algebras},
J. Algebra  {\bf 276} (2004), 776--793.

\bibitem{w2} J. ~Weyman, 
{\em Cohomology of vector bundles and syzygies},
London Mathematical Society Lecture Note Series, Cambridge, 2003.
 
\end{thebibliography}

\end{document}